\documentstyle[11pt]{article}
\begin{document}
\newcommand{\p}{\parallel }
\makeatletter \makeatother
\newtheorem{th}{Theorem}[section]
\newtheorem{lem}{Lemma}[section]
\newtheorem{de}{Definition}[section]
\newtheorem{rem}{Remark}[section]
\newtheorem{cor}{Corollary}[section]
\renewcommand{\theequation}{\thesection.\arabic {equation}}

\title{{\bf A general type of twisted anomaly cancellation formulas}}

\author{ Yong Wang\\
 }

\date{}
\maketitle

\begin{abstract}~~ For even dimensional manifolds, we prove some twisted anomaly cancellation formulas
which generalize some well-known cancellation formulas. For odd
dimensional manifolds, we obtain some modularly invariant
characteristic forms by the Chern-Simons transgression and we also
get some twisted anomaly cancellation
formulas.\\

 \noindent{\bf Subj. Class.:}\quad Differential geometry; Algebraic topology
\end{abstract}

\section{Introduction}
  \quad In 1983, the physicists Alvarez-Gaum\'{e} and Witten [AW]
  discovered the "miraculous cancellation" formula for gravitational
  anomaly which reveals a beautiful relation between the top
  components of the Hirzebruch $\widehat{L}$-form and
  $\widehat{A}$-form of a $12$-dimensional smooth Riemannian
  manifold. Kefeng Liu [Li] established higher dimensional "miraculous cancellation"
  formulas for $(8k+4)$-dimensional Riemannian manifolds by
  developing modular invariance properties of characteristic forms.
  These formulas could be used to deduce some divisibility results. In
  [HZ1], [HZ2], for each $(8k+4)$-dimensional smooth Riemannian
  manifold, a more general cancellation formula that involves a
  complex line bundle was established. This formula was applied to
  ${\rm spin}^c$ manifolds, then an analytic Ochanine congruence
  formula was derived. For $(8k+2)$ and $(8k+6)$-dimensional smooth Riemannian
  manifolds, F. Han and X. Huang [HH] obtained some cancellation formulas. They also got a general
  type of cancellation formulas.  \\
  \indent On the other hand, motivated by the Chern-Simons theory, in
  [CH], Qingtao Chen and Fei Han computed the transgressed forms of some modularly invariant characteristic forms,
which are related to the elliptic genera. They studied the
modularity properties of these secondary characteristic forms and
relations among them.  They also got an anomaly cancellation formula
for $11$-dimensional manifolds. In [W], the author computed the
transgressed forms of some modularly invariant characteristic forms,
which are related to the "twisted" elliptic genera and studied the
modularity properties of these secondary characteristic forms and
relations among them.  We also got some twisted anomaly cancellation
formulas on some odd dimensional manifolds. The purpose of paper is
to prove more general cancellation formulas for even and odd
dimensional manifolds. We hope that these new general cancellation
formulas obtained here
could be applied somewhere.\\
\indent This paper is organized as follows: In Section 2, we review
some knowledge on characteristic forms and modular forms that we are
going to use. In Section 3, we prove some general cancellation
formulas which involve two complex line bundles and generalize some
well-known cancellation formulas for even dimensional manifolds. In
Section 4, we apply the Chern-Simons transgression to characteristic
forms with modularity properties which are related to the "twisted"
elliptic genera and obtain some interesting secondary characteristic
forms with modularity properties. We also get two twisted
cancellation
formulas for $9$- and $11$-dimensional manifolds.\\

\section{characteristic forms and modular forms}
 \quad The purpose of this section is to review the necessary knowledge on
characteristic forms and modular forms that we are going to use.\\

 \noindent {\bf  2.1 characteristic forms. }Let $M$ be a Riemannian manifold.
 Let $\nabla^{ TM}$ be the associated Levi-Civita connection on $TM$
 and $R^{TM}=(\nabla^{TM})^2$ be the curvature of $\nabla^{ TM}$.
 Let $\widehat{A}(TM,\nabla^{ TM})$ and $\widehat{L}(TM,\nabla^{ TM})$
 be the Hirzebruch characteristic forms defined respectively by (cf.
 [Z])
 $$\widehat{A}(TM,\nabla^{ TM})={\rm
 det}^{\frac{1}{2}}\left(\frac{\frac{\sqrt{-1}}{4\pi}R^{TM}}{{\rm
 sinh}(\frac{\sqrt{-1}}{4\pi}R^{TM})}\right),$$
 $$\widehat{L}(TM,\nabla^{ TM})={\rm
 det}^{\frac{1}{2}}\left(\frac{\frac{\sqrt{-1}}{2\pi}R^{TM}}{{\rm
 tanh}(\frac{\sqrt{-1}}{4\pi}R^{TM})}\right).\eqno(2.1)$$
   Let $E$, $F$ be two Hermitian vector bundles over $M$ carrying
   Hermitian connection $\nabla^E,\nabla^F$ respectively. Let
   $R^E=(\nabla^E)^2$ (resp. $R^F=(\nabla^F)^2$) be the curvature of
   $\nabla^E$ (resp. $\nabla^F$). If we set the formal difference
   $G=E-F$, then $G$ carries an induced Hermitian connection
   $\nabla^G$ in an obvious sense. We define the associated Chern
   character form as
   $${\rm ch}(G,\nabla^G)={\rm tr}\left[{\rm
   exp}(\frac{\sqrt{-1}}{2\pi}R^E)\right]-{\rm tr}\left[{\rm
   exp}(\frac{\sqrt{-1}}{2\pi}R^F)\right].\eqno(2.2)$$
   For any complex number $t$, let
   $$\wedge_t(E)={\bf C}|_M+tE+t^2\wedge^2(E)+\cdots,~S_t(E)={\bf
   C}|_M+tE+t^2S^2(E)+\cdots$$
   denote respectively the total exterior and symmetric powers of
   $E$, which live in $K(M)[[t]].$ The following relations between
   these operations hold,
   $$S_t(E)=\frac{1}{\wedge_{-t}(E)},~\wedge_t(E-F)=\frac{\wedge_t(E)}{\wedge_t(F)}.\eqno(2.3)$$
   Moreover, if $\{\omega_i\},\{\omega_j'\}$ are formal Chern roots
   for Hermitian vector bundles $E,F$ respectively, then
   $${\rm ch}(\wedge_t(E))=\prod_i(1+e^{\omega_i}t).\eqno(2.4)$$
   Then we have the following formulas for Chern character forms,
   $${\rm ch}(S_t(E))=\frac{1}{\prod_i(1-e^{\omega_i}t)},~
{\rm
ch}(\wedge_t(E-F))=\frac{\prod_i(1+e^{\omega_i}t)}{\prod_j(1+e^{\omega_j'}t)}.\eqno(2.5)$$
\indent If $W$ is a real Euclidean vector bundle over $M$ carrying a
Euclidean connection $\nabla^W$, then its complexification $W_{\bf
C}=W\otimes {\bf C}$ is a complex vector bundle over $M$ carrying a
canonical induced Hermitian metric from that of $W$, as well as a
Hermitian connection $\nabla^{W_{\bf C}}$ induced from $\nabla^W$.
If $E$ is a vector bundle (complex or real) over $M$, set
$\widetilde{E}=E-{\rm dim}E$ in $K(M)$ or $KO(M)$.\\

\noindent{\bf 2.2 Some properties about the Jacobi theta functions
and modular forms}\\
   \indent We first recall the four Jacobi theta functions are
   defined as follows( cf. [Ch]):
   $$\theta(v,\tau)=2q^{\frac{1}{8}}{\rm sin}(\pi
   v)\prod_{j=1}^{\infty}[(1-q^j)(1-e^{2\pi\sqrt{-1}v}q^j)(1-e^{-2\pi\sqrt{-1}v}q^j)],\eqno(2.6)$$
$$\theta_1(v,\tau)=2q^{\frac{1}{8}}{\rm cos}(\pi
   v)\prod_{j=1}^{\infty}[(1-q^j)(1+e^{2\pi\sqrt{-1}v}q^j)(1+e^{-2\pi\sqrt{-1}v}q^j)],\eqno(2.7)$$
$$\theta_2(v,\tau)=\prod_{j=1}^{\infty}[(1-q^j)(1-e^{2\pi\sqrt{-1}v}q^{j-\frac{1}{2}})
(1-e^{-2\pi\sqrt{-1}v}q^{j-\frac{1}{2}})],\eqno(2.8)$$
$$\theta_3(v,\tau)=\prod_{j=1}^{\infty}[(1-q^j)(1+e^{2\pi\sqrt{-1}v}q^{j-\frac{1}{2}})
(1+e^{-2\pi\sqrt{-1}v}q^{j-\frac{1}{2}})],\eqno(2.9)$$ \noindent
where $q=e^{2\pi\sqrt{-1}\tau}$ with $\tau\in\textbf{H}$, the upper
half complex plane. Let
$$\theta'(0,\tau)=\frac{\partial\theta(v,\tau)}{\partial
v}|_{v=0}.\eqno(2.10)$$ \noindent Then the following Jacobi identity
(cf. [Ch]) holds,
$$\theta'(0,\tau)=\pi\theta_1(0,\tau)\theta_2(0,\tau)\theta_3(0,\tau).\eqno(2.11)$$
\noindent Denote $SL_2({\bf Z})=\left\{\left(\begin{array}{cc}
\ a & b  \\
 c  & d
\end{array}\right)\mid a,b,c,d \in {\bf Z},~ad-bc=1\right\}$ the
modular group. Let $S=\left(\begin{array}{cc}
\ 0 & -1  \\
 1  & 0
\end{array}\right),~T=\left(\begin{array}{cc}
\ 1 &  1 \\
 0  & 1
\end{array}\right)$ be the two generators of $SL_2(\bf{Z})$. They
act on $\textbf{H}$ by $S\tau=-\frac{1}{\tau},~T\tau=\tau+1$. One
has the following transformation laws of theta functions under the
actions of $S$ and $T$ (cf. [Ch]):
$$\theta(v,\tau+1)=e^{\frac{\pi\sqrt{-1}}{4}}\theta(v,\tau),~~\theta(v,-\frac{1}{\tau})
=\frac{1}{\sqrt{-1}}\left(\frac{\tau}{\sqrt{-1}}\right)^{\frac{1}{2}}e^{\pi\sqrt{-1}\tau
v^2}\theta(\tau v,\tau);\eqno(2.12)$$
$$\theta_1(v,\tau+1)=e^{\frac{\pi\sqrt{-1}}{4}}\theta_1(v,\tau),~~\theta_1(v,-\frac{1}{\tau})
=\left(\frac{\tau}{\sqrt{-1}}\right)^{\frac{1}{2}}e^{\pi\sqrt{-1}\tau
v^2}\theta_2(\tau v,\tau);\eqno(2.13)$$
$$\theta_2(v,\tau+1)=\theta_3(v,\tau),~~\theta_2(v,-\frac{1}{\tau})
=\left(\frac{\tau}{\sqrt{-1}}\right)^{\frac{1}{2}}e^{\pi\sqrt{-1}\tau
v^2}\theta_1(\tau v,\tau);\eqno(2.14)$$
$$\theta_3(v,\tau+1)=\theta_2(v,\tau),~~\theta_3(v,-\frac{1}{\tau})
=\left(\frac{\tau}{\sqrt{-1}}\right)^{\frac{1}{2}}e^{\pi\sqrt{-1}\tau
v^2}\theta_3(\tau v,\tau).\eqno(2.15)$$ \noindent Differentiating
the above transformation formulas, we get that
$$\theta'(v,\tau+1)=e^{\frac{\pi\sqrt{-1}}{4}}\theta'(v,\tau),$$
$$\theta'(v,-\frac{1}{\tau})=\frac{1}{\sqrt{-1}}\left(\frac{\tau}{\sqrt{-1}}\right)^{\frac{1}{2}}e^{\pi\sqrt{-1}\tau
v^2}(2\pi\sqrt{-1}\tau v\theta(\tau v,\tau)+\tau\theta'(\tau
v,\tau));$$
$$\theta'_1(v,\tau+1)=e^{\frac{\pi\sqrt{-1}}{4}}\theta_1'(v,\tau),$$
$$\theta'_1(v,-\frac{1}{\tau})=\left(\frac{\tau}{\sqrt{-1}}\right)^{\frac{1}{2}}e^{\pi\sqrt{-1}\tau
v^2}(2\pi\sqrt{-1}\tau v\theta_2(\tau v,\tau)+\tau\theta'_2(\tau
v,\tau));$$
$$\theta'_2(v,\tau+1)=\theta_3'(v,\tau),$$
$$\theta'_2(v,-\frac{1}{\tau})=\left(\frac{\tau}{\sqrt{-1}}\right)^{\frac{1}{2}}e^{\pi\sqrt{-1}\tau
v^2}(2\pi\sqrt{-1}\tau v\theta_1(\tau v,\tau)+\tau\theta'_1(\tau
v,\tau));$$
$$\theta'_3(v,\tau+1)=\theta_2'(v,\tau),$$
$$\theta'_3(v,-\frac{1}{\tau})=\left(\frac{\tau}{\sqrt{-1}}\right)^{\frac{1}{2}}e^{\pi\sqrt{-1}\tau
v^2}(2\pi\sqrt{-1}\tau v\theta_3(\tau v,\tau)+\tau\theta'_3(\tau
v,\tau))\eqno(2.16)$$
 \noindent Therefore
 $$\theta'(0,-\frac{1}{\tau})=\frac{1}{\sqrt{-1}}\left(\frac{\tau}{\sqrt{-1}}\right)^{\frac{1}{2}}
\tau\theta'(0,\tau).\eqno(2.17)$$
 \noindent {\bf Definition 2.1} A modular form over $\Gamma$, a
 subgroup of $SL_2({\bf Z})$, is a holomorphic function $f(\tau)$ on
 $\textbf{H}$ such that
 $$f(g\tau):=f\left(\frac{a\tau+b}{c\tau+d}\right)=\chi(g)(c\tau+d)^kf(\tau),
 ~~\forall g=\left(\begin{array}{cc}
\ a & b  \\
 c & d
\end{array}\right)\in\Gamma,\eqno(2.18)$$
\noindent where $\chi:\Gamma\rightarrow {\bf C}^{\star}$ is a
character of $\Gamma$. $k$ is called the weight of $f$.\\
Let $$\Gamma_0(2)=\left\{\left(\begin{array}{cc}
\ a & b  \\
 c  & d
\end{array}\right)\in SL_2({\bf Z})\mid c\equiv 0~({\rm
mod}~2)\right\},$$
$$\Gamma^0(2)=\left\{\left(\begin{array}{cc}
\ a & b  \\
 c  & d
\end{array}\right)\in SL_2({\bf Z})\mid b\equiv 0~({\rm
mod}~2)\right\},$$
$$\Gamma_\theta=\left\{\left(\begin{array}{cc}
\ a & b  \\
 c  & d
\end{array}\right)\in SL_2({\bf Z})\mid
\left(\begin{array}{cc}
\ a & b  \\
 c  & d
\end{array}\right)
\equiv \left(\begin{array}{cc}
\ 1 & 0  \\
 0  & 1
\end{array}\right)
{\rm or} \left(\begin{array}{cc}
\ 0 & 1  \\
 1  & 0
\end{array}\right)
 ~({\rm
mod}~2)\right\}$$ be the three modular subgroups of $SL_2({\bf Z})$.
It is known that the generators of $\Gamma_0(2)$ are $T,~ST^2ST$,
the generators of $\Gamma^0(2)$ are $STS,~T^2STS$ and the generators
of $\Gamma_\theta$ are $S,~T^2$ (cf.[Ch]).\\
\indent If $\Gamma$ is a modular subgroup, let ${\mathcal{M}}_{{\bf
R}}(\Gamma)$ denote the ring of modular forms over $\Gamma$ with
real Fourier coefficients. Writing $\theta_j=\theta_j(0,\tau),~1\leq
j\leq 3,$ we introduce six explicit modular forms (cf. [Li]),
$$\delta_1(\tau)=\frac{1}{8}(\theta_2^4+\theta_3^4),~~\varepsilon_1(\tau)=\frac{1}{16}\theta_2^4\theta_3^4,$$
$$\delta_2(\tau)=-\frac{1}{8}(\theta_1^4+\theta_3^4),~~\varepsilon_2(\tau)=\frac{1}{16}\theta_1^4\theta_3^4,$$
$$\delta_3(\tau)=\frac{1}{8}(\theta_1^4-\theta_2^4),~~\varepsilon_3(\tau)=-\frac{1}{16}\theta_1^4\theta_2^4.$$
\noindent They have the following Fourier expansions in
$q^{\frac{1}{2}}$:
$$\delta_1(\tau)=\frac{1}{4}+6q+\cdots,~~\varepsilon_1(\tau)=\frac{1}{16}-q+\cdots,$$
$$\delta_2(\tau)=-\frac{1}{8}-3q^{\frac{1}{2}}+\cdots,~~\varepsilon_2(\tau)=q^{\frac{1}{2}}+\cdots,$$
$$\delta_3(\tau)=-\frac{1}{8}+3q^{\frac{1}{2}}+\cdots,~~\varepsilon_3(\tau)=-q^{\frac{1}{2}}+\cdots,$$
\noindent where the $"\cdots"$ terms are the higher degree terms,
all of which have integral coefficients. They also satisfy the
transformation laws,
$$\delta_2(-\frac{1}{\tau})=\tau^2\delta_1(\tau),~~~~~~\varepsilon_2(-\frac{1}{\tau})
=\tau^4\varepsilon_1(\tau),\eqno(2.19)$$
$$\delta_2(\tau+1)=\delta_3(\tau),~~~~~~\varepsilon_2(\tau+1)=\varepsilon_3(\tau).\eqno(2.20)$$
\noindent {\bf Lemma 2.2} ([Li]) {\it $\delta_1(\tau)$ (resp.
$\varepsilon_1(\tau)$) is a modular form of weight $2$ (resp. $4$)
over $\Gamma_0(2)$, $\delta_2(\tau)$ (resp. $\varepsilon_2(\tau)$)
is a modular form of weight $2$ (resp. $4$) over $\Gamma^0(2)$,
while  $\delta_3(\tau)$ (resp. $\varepsilon_3(\tau)$) is a modular
form of weight $2$ (resp. $4$) over $\Gamma_\theta(2)$ and moreover
${\mathcal{M}}_{{\bf R}}(\Gamma^0(2))={\bf
R}[\delta_2(\tau),\varepsilon_2(\tau)]$.}

\section {A general type of cancellation formulas for even dimensional manifolds}

   \quad Let $M$ be a $2d$ dimensional Riemannian
manifold and
 $\xi^0,$~ $\xi$ be rank two real oriented Euclidean vector
   bundles over $M$ carrying with Euclidean connections
   $\nabla^{\xi^0}$, $\nabla^\xi$. Set
   $$\Theta_1(T_{C}M,m_0\xi^0_C,\xi_C)=
   \bigotimes _{n=1}^{\infty}S_{q^n}(\widetilde{T_CM}-m_0\widetilde{\xi^0_C})\otimes
\bigotimes
_{m=1}^{\infty}\wedge_{q^m}(\widetilde{T_CM}-m_0\widetilde{\xi^0_C}-2\widetilde{\xi_C})$$
$$~~~~~~~~\otimes \bigotimes _{r=1}^{\infty}\wedge
_{q^{r-\frac{1}{2}}}(\widetilde{\xi_C})\otimes\bigotimes
_{s=1}^{\infty}\wedge _{-q^{s-\frac{1}{2}}}(\widetilde{\xi_C}),$$
$$\Theta_2(T_{C}M,m_0\xi^0_C,\xi_C)=\bigotimes _{n=1}^{\infty}S_{q^n}(\widetilde{T_CM}-m_0\widetilde{\xi^0_C})\otimes
\bigotimes
_{m=1}^{\infty}\wedge_{-q^{m-\frac{1}{2}}}(\widetilde{T_CM}-m_0\widetilde{\xi^0_C}-2\widetilde{\xi_C})$$
$$~~~~~~~~\otimes \bigotimes _{r=1}^{\infty}\wedge
_{q^{r-\frac{1}{2}}}(\widetilde{\xi_C})\otimes\bigotimes
_{s=1}^{\infty}\wedge _{q^{s}}(\widetilde{\xi_C}),\eqno(3.1)$$
Clearly, $\Theta_1(T_{C}M,m_0\xi^0_C,\xi_C)$ and
$\Theta_2(T_{C}M,m_0\xi^0_C,\xi_C)$ admit formal Fourier expansion
in $q^{\frac{1}{2}}$ as
$$\Theta_1(T_{C}M,m_0\xi^0_C,\xi_C)=A_0(T_{C}M,m_0\xi^0_C,\xi_C)+A_1(T_{C}M,m_0\xi^0_C,\xi_C)q
^{\frac{1}{2}}+\cdots,$$
$$\Theta_2(T_{C}M,m_0\xi^0_C,\xi_C)=B_0(T_{C}M,m_0\xi^0_C,\xi_C)+B_1(T_{C}M,m_0\xi^0_C,\xi_C)q
^{\frac{1}{2}}+\cdots,\eqno(3.2)$$ where the $A_j$ and $B_j$ are
elements in the semi-group formally generated by Hermitian vector
bundles over $M$. Moreover, they carry canonically induced Hermitian
connections.
 Let $c_0=e(\xi,\nabla^{\xi^0})$ and
$c=e(\xi,\nabla^{\xi})$ be the Euler forms of $\xi^0$, $\xi$
canonically associated to $\nabla^{\xi^0}$~$\nabla^\xi$
respectively. If $\omega$ is a differential form over $M$, we denote
$\omega^{(2d)}$ its top degree component. Let $n$ be a nonnegative
integer and satisfy $d-\left(2n+\frac{1-(-1)^d}{2}\right)>0$,
then we have\\

\noindent {\bf Theorem 3.1} {\it The following identity holds,}
$$\left\{\frac{\widehat{L}(TM,\nabla^{TM})}{{\rm
cosh}^2{\frac{c}{2}}}\frac{({\rm
sinh}\frac{c_0}{2})^{2n+\frac{1-(-1)^d}{2}}}{({\rm
cosh}\frac{c_0}{2})^{2n+\frac{1-(-1)^d}{2}}}\right\}^{(2d)}$$
$$=
2^{\frac{3}{2}d-n-\frac{1-(-1)^d}{4}}\sum_{r=0}^{[\frac{m_1}{2}]}2^{-6r}\left\{
\widehat{A}(TM,\nabla^{TM}){\rm cosh}{\frac{c}{2}}({\rm
sinh}\frac{c_0}{2})^{2n+\frac{1-(-1)^d}{2}}\right.$$
$$\left.\cdot{\rm
ch}(b_r(T_CM,(2n+\frac{1-(-1)^d}{2})\xi^0_C,\xi_C))\right\}^{(2d)}
,\eqno(3.3)$$ {\it where $m_1=\frac{d}{2}-n-\frac{1-(-1)^d}{4}$ and
each $b_r(T_CM,(2n+\frac{1-(-1)^d}{2})\xi^0_C,\xi_C)$,~$0\leq r\leq
[\frac{m_1}{2}]$, is a canonical integral linear combination of
$B_j(T_CM,(2n+\frac{1-(-1)^d}{2})\xi^0_C,\xi_C)$, $0\leq j\leq
r.$}\\

\noindent {\bf Proof.} Let $\{\pm2\pi\sqrt{-1}x_j| ~1\leq j\leq d\}$
 be the Chern roots of $T_CM$ and $c_0=2\pi\sqrt{-1}u'.$
$c=2\pi\sqrt{-1}u.$ Set
$$Q_1(\tau)=\frac{\widehat{L}(TM,\nabla^{TM})}{{\rm
cosh}^2{\frac{c}{2}}}\frac{({\rm
sinh}\frac{c_0}{2})^{2n+\frac{1-(-1)^d}{2}}}{({\rm
cosh}\frac{c_0}{2})^{2n+\frac{1-(-1)^d}{2}}} {\rm
ch}(\Theta_1(T_CM,(2n+\frac{1-(-1)^d}{2})\xi^0_C,\xi_C)),\eqno(3.4)$$
$$Q_2(\tau)=\widehat{A}(TM,\nabla^{TM}){\rm cosh}{\frac{c}{2}}({\rm
sinh}\frac{c_0}{2})^{2n+\frac{1-(-1)^d}{2}}{\rm
ch}(\Theta_2(T_CM,(2n+\frac{1-(-1)^d}{2})\xi^0_C,\xi_C)),\eqno(3.5)$$
Let $\Theta_1(T_{C}M,\xi_C)=\Theta_1(T_{C}M,m_0C^2,\xi_C),~
\Theta_2(T_{C}M,\xi_C)=\Theta_2(T_{C}M,m_0C^2,\xi_C).$ Then
$$Q_1(\tau)=\frac{\widehat{L}(TM,\nabla^{TM})}{{\rm
cosh}^2{\frac{c}{2}}}{\rm ch}(\Theta_1(T_CM,,\xi_C))$$
$$\cdot\left[\frac{{\rm cosh}{\frac{c_0}{2}}}{{\rm
sinh}\frac{c_0}{2}}{\rm ch} \left(\bigotimes
_{n=1}^{\infty}S_{q^n}(\widetilde{\xi^0_C})\otimes \bigotimes
_{m=1}^{\infty}\wedge_{q^m}(\widetilde{\xi^0_C})\right)\right]^{-2n-\frac{1-(-1)^d}{2}}.\eqno(3.6)$$
By Proposition 2.5 in [HZ2], we have
$$\frac{\widehat{L}(TM,\nabla^{TM})}{{\rm
cosh}^2{\frac{c}{2}}}{\rm ch}(\Theta_1(T_CM,,\xi_C))~~~~~~~~~~~~~~$$
$$ =2^d \left\{
\prod_{j=1}^{d}\left(x_j\frac{\theta'(0,\tau)}{\theta(x_j,\tau)}\frac
{\theta_1(x_j,\tau)}{\theta_1(0,\tau)}\right)\frac{\theta_1^2(0,\tau)}
{\theta_1^2(u,\tau)}\frac{\theta_3(u,\tau)}{\theta_3(0,\tau)}
\frac{\theta_2(u,\tau)}{\theta_2(0,\tau)}\right\}.\eqno(3.7)$$
Direct computations show that $$ \frac{{\rm
cosh}{\frac{c_0}{2}}}{{\rm sinh}\frac{c_0}{2}}{\rm ch}
\left(\bigotimes _{n=1}^{\infty}S_{q^n}(\widetilde{\xi^0_C})\otimes
\bigotimes _{m=1}^{\infty}\wedge_{q^m}(\widetilde{\xi^0_C})\right)
=\frac{1}{\pi\sqrt{-1}}\frac{\theta'(0,\tau)}{\theta(u',\tau)}\frac
{\theta_1(u',\tau)}{\theta_1(0,\tau)}.\eqno(3.8)$$ By (3.6)-(3.8),
we have
$$ Q_1(\tau)=2^d (\pi\sqrt{-1})^{2n+\frac{1-(-1)^d}{2}}\left\{
\prod_{j=1}^{d}\left(x_j\frac{\theta'(0,\tau)}{\theta(x_j,\tau)}\frac
{\theta_1(x_j,\tau)}{\theta_1(0,\tau)}\right)\right.$$
$$\left.\left(\frac{\theta(u',\tau)}{\theta'(0,\tau)}\frac
{\theta_1(0,\tau)}{\theta_1(u',\tau)}\right)^{2n+\frac{1-(-1)^d}{2}}
\cdot \frac{\theta_1^2(0,\tau)}
{\theta_1^2(u,\tau)}\frac{\theta_3(u,\tau)}{\theta_3(0,\tau)}
\frac{\theta_2(u,\tau)}{\theta_2(0,\tau)}\right\}.\eqno(3.9)$$
Similarly,
$$Q_2(\tau)=\widehat{A}(TM,\nabla^{TM}){\rm
cosh}{\frac{c}{2}}{\rm ch}(\Theta_2(T_CM,\xi_C))$$ $$\cdot\left[
({\rm sinh}\frac{c_0}{2}) {\rm ch} \left(\bigotimes
_{n=1}^{\infty}S_{q^n}(-\widetilde{\xi^0_C})\otimes \bigotimes
_{m=1}^{\infty}\wedge_{-q^{m-\frac{1}{2}}}(-\widetilde{\xi^0_C})\right)\right]^{2n+\frac{1-(-1)^d}{2}};
\eqno(3.10)$$
$$\widehat{A}(TM,\nabla^{TM}){\rm cosh}{\frac{c}{2}}{\rm
ch}(\Theta_2(T_CM,\xi_C))$$
$$=\left(\prod_{j=1}^{d}x_j\frac{\theta'(0,\tau)}{\theta(x_j,\tau)}\frac
{\theta_2(x_j,\tau)}{\theta_2(0,\tau)}\right)\frac{\theta_2^2(0,\tau)}
{\theta_2^2(u,\tau)}\frac{\theta_3(u,\tau)}{\theta_3(0,\tau)}
\frac{\theta_1(u,\tau)}{\theta_1(0,\tau)};\eqno(3.11)$$
$$
{\rm sinh}\frac{c_0}{2} {\rm ch} \left(\bigotimes
_{n=1}^{\infty}S_{q^n}(-\widetilde{\xi^0_C})\otimes \bigotimes
_{m=1}^{\infty}\wedge_{-q^{m-\frac{1}{2}}}(-\widetilde{\xi^0_C})\right)
=\sqrt{-1}\pi\frac{\theta(u',\tau)}{\theta'(0,\tau)}\frac
{\theta_2(0,\tau)}{\theta_2(u',\tau)},\eqno(3.12)$$ so we have
$$Q_2(\tau)=(\sqrt{-1}\pi)^{2n+\frac{1-(-1)^d}{2}}
\left(\prod_{j=1}^{d}x_j\frac{\theta'(0,\tau)}{\theta(x_j,\tau)}\frac
{\theta_2(x_j,\tau)}{\theta_2(0,\tau)}\right)$$
$$\cdot\left(\frac{\theta(u',\tau)}{\theta'(0,\tau)}\frac
{\theta_2(0,\tau)}{\theta_2(u',\tau)}\right)^{2n+\frac{1-(-1)^d}{2}}
\frac{\theta_2^2(0,\tau)}
{\theta_2^2(u,\tau)}\frac{\theta_3(u,\tau)}{\theta_3(0,\tau)}
\frac{\theta_1(u,\tau)}{\theta_1(0,\tau)}.\eqno(3.13)$$ Let
$P_1(\tau)=Q_1(\tau)^{(2d)},~P_2(\tau)=Q_2(\tau)^{(2d)}.$ By
(2.12)-(2.15) and (2.17), then $P_1(\tau)$ is a modular form of
weight $d-(2n+\frac{1-(-1)^d}{2})$ over $\Gamma_0(2)$, while
$P_2(\tau)$ is a modular form of weight $d-(2n+\frac{1-(-1)^d}{2})$
over $\Gamma^0(2)$ . Moreover, the following identity holds,
$$P_1(-\frac{1}{\tau})=2^d\tau^{d-(2n+\frac{1-(-1)^d}{2})}P_2(\tau).\eqno(3.14)$$
\indent Observe that at any point $x\in M$, up to the volume form
determined by the metric on $T_xM$, both $P_i(\tau),$ $i=1,2$, can
be view as a power series of $q^{\frac{1}{2}}$ with real Fourier
coefficients. By Lemma 2.2, we have
$$P_2(\tau)=h_0(8\delta_2)^{m_1}+h_1(8\delta_2)^{m_1-2}\varepsilon_2+\cdots+h_{[\frac{m_1}{2}]}(8\delta_2)
^{m_1-2[\frac{m_1}{2}]}\varepsilon^{[\frac{m_1}{2}]}_2,\eqno(3.15)$$
where each $h_j,~ 0\leq j\leq [\frac{m_1}{2}],$ is a real multiple
of the volume form at $x$. By (2.19) (3.14) and (3.15), we get
$$P_1(\tau)=2^d\left[h_0(8\delta_1)^{m_1}+h_1(8\delta_1)^{m_1-2}\varepsilon_1+\cdots+h_{[\frac{m_1}{2}]}(8\delta_1)
^{m_1-2[\frac{m_1}{2}]}\varepsilon^{[\frac{m_1}{2}]}_1\right].\eqno(3.16)$$
By comparing the constant term in (3.16), we get
$$\left\{\frac{\widehat{L}(TM,\nabla^{TM})}{{\rm
cosh}^2{\frac{c}{2}}}\frac{({\rm
sinh}\frac{c_0}{2})^{2n+\frac{1-(-1)^d}{2}}}{({\rm
cosh}\frac{c_0}{2})^{2n+\frac{1-(-1)^d}{2}}}\right\}^{(2d)} =
2^{\frac{3}{2}d-n-\frac{1-(-1)^d}{4}}\sum_{r=0}^{[\frac{m_1}{2}]}2^{-6r}h_r.\eqno(3.17)$$
By comparing the coefficients of $q^{\frac{j}{2}}$,~$j\geq 0$
between the two sides of (3.15), we can use the induction method to
prove that each $h_r~0\leq r\leq [\frac{m_1}{2}]$, can be expressed
through a canonical integral linear combination of
$$\left\{\widehat{A}(TM,\nabla^{TM}){\rm cosh}{\frac{c}{2}}({\rm
sinh}\frac{c_0}{2})^{2n+\frac{1-(-1)^d}{2}} {\rm
ch}(B_r(T_CM,(2n+\frac{1-(-1)^d}{2})\xi^0_C,\xi_C))\right\}^{(2d)}.$$
Here we write out the explicit expressions for $h_0$ and $h_1$ as
follows.
$$h_0=(-1)^{m_1}\left\{\widehat{A}(TM,\nabla^{TM}){\rm cosh}{\frac{c}{2}}({\rm
sinh}\frac{c_0}{2})^{2n+\frac{1-(-1)^d}{2}}\right\}^{(2d)},\eqno(3.18)$$
$$h_1=(-1)^{m_1}\left\{\widehat{A}(TM,\nabla^{TM}){\rm cosh}{\frac{c}{2}}({\rm
sinh}\frac{c_0}{2})^{2n+\frac{1-(-1)^d}{2}}\right.$$
$$\left.\cdot\left( {\rm
ch}(B_1(T_CM,(2n+\frac{1-(-1)^d}{2})\xi^0_C,\xi_C))-24m_1\right)\right\}.\eqno(3.19)$$
$\Box$ \\
\indent Putting $d=4k+2$ and $n=0$ in Theorem 3.1, we get
the
Han-Zhang cancellation formula (cf. [HZ2]),\\

 \noindent {\bf Corollary 3.2} {\it The following cancellation
 formula holds}
$$\left\{\frac{\widehat{L}(TM,\nabla^{TM})}{{\rm
cosh}^2{\frac{c}{2}}}\right\}^{(8k+4)}=8
\sum_{r=0}^{k}2^{6k-6r}\left\{ \widehat{A}(TM,\nabla^{TM}){\rm
cosh}{\frac{c}{2}}{\rm ch}(b_r(T_CM,\xi_C))\right\}^{(8k+4)}
.\eqno(3.20)$$

\indent If $\xi$ is a trivial bundle, we get the
Han-Huang cancellation formula (cf. [HH]),\\

\noindent {\bf Corollary 3.3} {\it The following cancellation
 formula holds}
$$\left\{\widehat{L}(TM,\nabla^{TM})\frac{({\rm
sinh}\frac{c_0}{2})^{2n+\frac{1-(-1)^d}{2}}}{({\rm
cosh}\frac{c_0}{2})^{2n+\frac{1-(-1)^d}{2}}}\right\}^{(2d)}$$
$$=
2^{\frac{3}{2}d-n-\frac{1-(-1)^d}{4}}\sum_{r=0}^{[\frac{m_1}{2}]}2^{-6r}\left\{
\widehat{A}(TM,\nabla^{TM})({\rm
sinh}\frac{c_0}{2})^{2n+\frac{1-(-1)^d}{2}}\right.$$
$$\left.\cdot{\rm
ch}(b_r(T_CM,(2n+\frac{1-(-1)^d}{2})\xi^0_C,C^2))\right\}^{(2d)}
,\eqno(3.21)$$\\

\indent Putting $d=6$ and $n=1$, i.e. for $12$-dimensional manifold
$M$, we
have\\

\noindent {\bf Corollary 3.4} {\it The following cancellation
 formula holds}
$$\left\{\frac{\widehat{L}(TM,\nabla^{TM})}{{\rm
cosh}^2{\frac{c}{2}}}\frac{({\rm sinh}\frac{c_0}{2})^{2}}{({\rm
cosh}\frac{c_0}{2})^{2}}\right\}^{(12)} = \left\{
\widehat{A}(TM,\nabla^{TM}){\rm cosh}{\frac{c}{2}}({\rm
sinh}\frac{c_0}{2})^{2}\right.$$
$$\left.\cdot\left(112-4{\rm
ch}(T_CM,\nabla^{T_CM})+8(e^{c_0}+e^{-c_0}-2)+12(e^{c}+e^{-c}-2)\right)\right\}^{(12)}.\eqno(3.22)$$
\\
\indent Putting $d=6$ and $n=2$, i.e. for $12$-dimensional manifold
$M$, we have\\

\noindent {\bf Corollary 3.5} {\it The following cancellation
 formula holds}
$$\left\{\frac{\widehat{L}(TM,\nabla^{TM})}{{\rm
cosh}^2{\frac{c}{2}}}\frac{({\rm sinh}\frac{c_0}{2})^{4}}{({\rm
cosh}\frac{c_0}{2})^{4}}\right\}^{(12)} =-128\left\{
\widehat{A}(TM,\nabla^{TM}){\rm cosh}{\frac{c}{2}}({\rm
sinh}\frac{c_0}{2})^{4}\right\}^{(12)}.\eqno(3.23)$$
\\
\indent Putting $d=5$ and $n=0$, i.e. for $10$-dimensional manifold
$M$, we have\\

\noindent {\bf Corollary 3.6} {\it The following cancellation
 formula holds}
$$\left\{\frac{\widehat{L}(TM,\nabla^{TM})}{{\rm
cosh}^2{\frac{c}{2}}}\frac{({\rm sinh}\frac{c_0}{2})}{({\rm
cosh}\frac{c_0}{2})}\right\}^{(10)} =\left\{
\widehat{A}(TM,\nabla^{TM}){\rm cosh}{\frac{c}{2}}({\rm
sinh}\frac{c_0}{2}) \right.$$
$$\left.\cdot\left(52-2{\rm
ch}(T_CM,\nabla^{T_CM})+2(e^{c_0}+e^{-c_0}-2)+6(e^{c}+e^{-c}-2)\right)\right\}^{(10)}.\eqno(3.24)$$
\\
\indent Putting $d=5$ and $n=1$, i.e. for $10$-dimensional manifold
$M$, we have\\

\noindent {\bf Corollary 3.6} {\it The following cancellation
 formula holds}
$$\left\{\frac{\widehat{L}(TM,\nabla^{TM})}{{\rm
cosh}^2{\frac{c}{2}}}\frac{({\rm sinh}\frac{c_0}{2})^{3}}{({\rm
cosh}\frac{c_0}{2})^{3}}\right\}^{(10)} =-64\left\{
\widehat{A}(TM,\nabla^{TM}){\rm cosh}{\frac{c}{2}}({\rm
sinh}\frac{c_0}{2})^{3}\right\}^{(10)}.\eqno(3.25)$$
\\
\indent Nextly we go on to prove some cancellation formulas. Define
$$\Theta_1(T_{C}M+\xi_C,m_0\xi^0_C,\xi_C)=
   \bigotimes _{n=1}^{\infty}S_{q^n}(\widetilde{T_CM}+\widetilde{\xi_C}-m_0\widetilde{\xi^0_C})$$ $$\otimes
\bigotimes
_{m=1}^{\infty}\wedge_{q^m}(\widetilde{T_CM}+\widetilde{\xi_C}-m_0\widetilde{\xi^0_C}-2\widetilde{\xi_C})
\otimes \bigotimes _{r=1}^{\infty}\wedge
_{q^{r-\frac{1}{2}}}(\widetilde{\xi_C})\otimes\bigotimes
_{s=1}^{\infty}\wedge _{-q^{s-\frac{1}{2}}}(\widetilde{\xi_C}),$$
$$\Theta_2(T_{C}M+\xi_C,m_0\xi^0_C,\xi_C)=\bigotimes _{n=1}^{\infty}S_{q^n}(\widetilde{T_CM}+\widetilde{\xi_C}-m_0
\widetilde{\xi^0_C})$$ $$\otimes \bigotimes
_{m=1}^{\infty}\wedge_{-q^{m-\frac{1}{2}}}(\widetilde{T_CM}+\widetilde{\xi_C}
-m_0\widetilde{\xi^0_C}-2\widetilde{\xi_C}) \otimes \bigotimes
_{r=1}^{\infty}\wedge
_{q^{r-\frac{1}{2}}}(\widetilde{\xi_C})\otimes\bigotimes
_{s=1}^{\infty}\wedge _{q^{s}}(\widetilde{\xi_C}),\eqno(3.26)$$
$\Theta_1(T_{C}M+\xi_C,m_0\xi^0_C,\xi_C)$ and
$\Theta_2(T_{C}M+\xi_C,m_0\xi^0_C,\xi_C)$ admit formal Fourier
expansion in $q^{\frac{1}{2}}$ as
$$\Theta_1(T_{C}M+\xi_C,m_0\xi^0_C,\xi_C)=A'_0(T_{C}M,m_0\xi^0_C,\xi_C)+A_1'(T_{C}M,m_0\xi^0_C,\xi_C)q
^{\frac{1}{2}}+\cdots,$$
$$\Theta_2(T_{C}M+\xi_C,m_0\xi^0_C,\xi_C)=B'_0(T_{C}M,m_0\xi^0_C,\xi_C)+B_1'(T_{C}M,m_0\xi^0_C,\xi_C)q
^{\frac{1}{2}}+\cdots,\eqno(3.27)$$ Set
$$Q'_1(\tau)=\widehat{L}(TM,\nabla^{TM})\frac{{\rm
cosh}{\frac{c}{2}}}{{\rm sinh}{\frac{c}{2}}}\frac{({\rm
sinh}\frac{c_0}{2})^{2n+\frac{1+(-1)^d}{2}}}{({\rm
cosh}\frac{c_0}{2})^{2n+\frac{1+(-1)^d}{2}}}$$
$$ \cdot\left({\rm
ch}(\Theta_1(T_CM+\xi_C,(2n+\frac{1+(-1)^d}{2})\xi^0_C,C^2))\right.$$
$$\left.-\frac{{\rm
ch}(\Theta_1(T_CM+\xi_C,(2n+\frac{1+(-1)^d}{2})\xi^0_C,\xi_C)}{{\rm
cosh}^2{\frac{c}{2}}}\right),\eqno(3.28)$$
$$Q'_2(\tau)=\widehat{A}(TM,\nabla^{TM})\frac{1}{2{\rm sinh}\frac{c}{2}}({\rm
sinh}\frac{c_0}{2})^{2n+\frac{1+(-1)^d}{2}}$$
$$\cdot\left({\rm
ch}(\Theta_2(T_CM+\xi_C,(2n+\frac{1+(-1)^d}{2})\xi^0_C,C^2))\right.$$
$$\left.-{\rm cosh}(\frac{c}{2}){\rm
ch}(\Theta_2(T_CM+\xi_C,(2n+\frac{1+(-1)^d}{2})\xi^0_C,\xi_C))\right).\eqno(3.29)$$
\noindent Direct computations show that
$$ Q'_1(\tau)=2^d (\pi\sqrt{-1})^{2n+\frac{1+(-1)^d}{2}-1}\left(
\prod_{j=1}^{d}x_j\frac{\theta'(0,\tau)}{\theta(x_j,\tau)}\frac
{\theta_1(x_j,\tau)}{\theta_1(0,\tau)}\right)
\left(\frac{\theta(u',\tau)}{\theta'(0,\tau)}\frac
{\theta_1(0,\tau)}{\theta_1(u',\tau)}\right)^{2n+\frac{1+(-1)^d}{2}}$$
$$\cdot\frac{\theta'(0,\tau)}{\theta(u,\tau)}
\left(\frac {\theta_1(u,\tau)}{\theta_1(0,\tau)} -\frac
{\theta_1(0,\tau)}{\theta_1(u,\tau)} \frac
{\theta_3(u,\tau)}{\theta_3(0,\tau)} \frac
{\theta_2(u,\tau)}{\theta_2(0,\tau)}\right),\eqno(3.30)$$
$$Q'_2(\tau)=\frac{1}{2}(\sqrt{-1}\pi)^{2n+\frac{1+(-1)^d}{2}-1}
\left(\prod_{j=1}^{d}x_j\frac{\theta'(0,\tau)}{\theta(x_j,\tau)}\frac
{\theta_2(x_j,\tau)}{\theta_2(0,\tau)}\right)
\left(\frac{\theta(u',\tau)}{\theta'(0,\tau)}\frac
{\theta_2(0,\tau)}{\theta_2(u',\tau)}\right)^{2n+\frac{1+(-1)^d}{2}}$$
$$\cdot\frac{\theta'(0,\tau)}{\theta(u,\tau)}\left(\frac {\theta_2(u,\tau)}{\theta_2(0,\tau)} -\frac
{\theta_2(0,\tau)}{\theta_2(u,\tau)} \frac
{\theta_3(u,\tau)}{\theta_3(0,\tau)} \frac
{\theta_1(u,\tau)}{\theta_1(0,\tau)}\right).\eqno(3.31)$$ Let
$P'_1(\tau)=Q'_1(\tau)^{(2d)},~P'_2(\tau)=Q'_2(\tau)^{(2d)}$,
similarly we have $P'_1(\tau)$ is a modular form of weight
$d+1-(2n+\frac{1+(-1)^d}{2})$ over $\Gamma_0(2)$, while $P'_2(\tau)$
is a modular form of weight $d+1-(2n+\frac{1+(-1)^d}{2})$ over
$\Gamma^0(2)$ . Moreover, the following identity holds,
$$P'_1(-\frac{1}{\tau})=2^{d+1}\tau^{d+1-(2n+\frac{1+(-1)^d}{2})}P_2'(\tau).\eqno(3.32)$$
Let $n$ is a nonnegative integer and satisfy
$d-1-(2n+\frac{1+(-1)^d}{2})>0$.  Using the same trick in the proof
of Theorem 3.1, we obtain\\

\noindent {\bf Theorem 3.7} {\it The following identity holds,}
$$\left\{\widehat{L}(TM,\nabla^{TM})\frac{{\rm
sinh}{\frac{c}{2}}}{{\rm cosh}{\frac{c}{2}}}\frac{({\rm
sinh}\frac{c_0}{2})^{2n+\frac{1+(-1)^d}{2}}}{({\rm
cosh}\frac{c_0}{2})^{2n+\frac{1+(-1)^d}{2}}}\right\}^{(2d)} =
2^{\frac{3}{2}(d+1)-n-\frac{1+(-1)^d}{4}}\sum_{r=0}^{[\frac{m_2}{2}]}2^{-6r}h_r,\eqno(3.33)$$
 {\it where $m_2=\frac{d+1}{2}-n-\frac{1+(-1)^d}{4}$
and each $h_r$,~$0\leq r\leq [\frac{m_2}{2}]$, is a canonical
integral linear combination of}
$$\left\{\widehat{A}(TM,\nabla^{TM})\frac{1}{2{\rm sinh}\frac{c}{2}}({\rm
sinh}\frac{c_0}{2})^{2n+\frac{1+(-1)^d}{2}} \left({\rm
ch}(B'_j(T_CM,(2n+\frac{1+(-1)^d}{2})\xi^0_C,C^2))\right.\right.$$
$$\left.\left.-{\rm cosh}(\frac{c}{2}){\rm
ch}(B_j'(T_CM,(2n+\frac{1+(-1)^d}{2})\xi^0_C,\xi_C))\right)\right\}^{(2d)},~0\leq
j\leq r.$$
\\
\indent Putting $d=4k+1$ and $n=0$ in Theorem 3.7, we get the
Han-Huang cancellation formula (cf. [HH]),\\

 \noindent {\bf Corollary 3.8} {\it The following cancellation
 formula holds}
$$\left\{\widehat{L}(TM,\nabla^{TM})\frac{{\rm
sinh}{\frac{c}{2}}}{{\rm cosh}{\frac{c}{2}}}\right\}^{(8k+2)} =
8\sum_{r=0}^{k}2^{6k-6r}h_r.\eqno(3.34)$$
\\
\indent Putting $d=6$ and $n=1$, i.e. for $12$-dimensional manifold
$M$, we have\\

\noindent {\bf Corollary 3.9} {\it The following cancellation
 formula holds}
$$\left\{\widehat{L}(TM,\nabla^{TM})\frac{{\rm
sinh}{\frac{c}{2}}}{{\rm cosh}{\frac{c}{2}}}\frac{({\rm
sinh}\frac{c_0}{2})^{3}}{({\rm
cosh}\frac{c_0}{2})^{3}}\right\}^{(12)}
=\left\{\widehat{A}(TM,\nabla^{TM})\frac{1}{2{\rm
sinh}\frac{c}{2}}({\rm sinh}\frac{c_0}{2})^{3}\right.$$
$$\cdot\left[\left(224+24(e^{c_0}+e^{-c_0}-2)-8{\rm
ch}(T_CM,\nabla^{T_CM})\right)(1-{\rm cosh}\frac{c}{2})\right.$$
$$\left.\left.-8(e^c+e^{-c}-2)(1+2{\rm
cosh}\frac{c}{2})\right]\right\}^{(12)}.\eqno(3.35)$$
\\
\indent Putting $d=5$ and $n=1$, i.e. for $10$-dimensional manifold
$M$, we have\\

\noindent {\bf Corollary 3.10} {\it The following cancellation
 formula holds}
$$\left\{\widehat{L}(TM,\nabla^{TM})\frac{{\rm
sinh}{\frac{c}{2}}}{{\rm cosh}{\frac{c}{2}}}\frac{({\rm
sinh}\frac{c_0}{2})^{2}}{({\rm
cosh}\frac{c_0}{2})^{2}}\right\}^{(10)}
=\left\{\widehat{A}(TM,\nabla^{TM})\frac{1}{2{\rm
sinh}\frac{c}{2}}({\rm sinh}\frac{c_0}{2})^{2}\right.$$
$$\cdot\left[\left(104+8(e^{c_0}+e^{-c_0}-2)-4{\rm
ch}(T_CM,\nabla^{T_CM})\right)(1-{\rm cosh}\frac{c}{2})\right.$$
$$\left.\left.-(e^c+e^{-c}-2)(1+2{\rm
cosh}\frac{c}{2})\right]\right\}^{(10)}.\eqno(3.36)$$
\\
\section{ Transgressed forms and modularities}

    \quad In this section, following [CH], we transgress the modular characteristic forms
    in Section 3 and then get some cancellation formulas.\\
    \indent Let $M$ be $(2d-1)$-dimensional manifold. Set\\
   $$\Theta_1(T_{C}M,m_0\xi^0_C)=\Theta_1(T_{C}M,m_0\xi^0_C,C^2);~
   \Theta_2(T_{C}M,m_0\xi^0_C)=\Theta_2(T_{C}M,m_0\xi^0_C,C^2); $$
$$\Theta_3(T_{C}M,m_0\xi^0_C)=
   \bigotimes _{n=1}^{\infty}S_{q^n}(\widetilde{T_CM}-m_0\widetilde{\xi^0_C})\otimes
\bigotimes
_{m=1}^{\infty}\wedge_{q^{m-\frac{1}{2}}}(\widetilde{T_CM}-m_0\widetilde{\xi^0_C}).$$
Set
$$\Phi_L(\tau)=\widehat{L}(TM,\nabla^{TM})\frac{({\rm
sinh}\frac{c_0}{2})^{2n+\frac{1-(-1)^d}{2}}}{({\rm
cosh}\frac{c_0}{2})^{2n+\frac{1-(-1)^d}{2}}} {\rm
ch}(\Theta_1(T_CM,(2n+\frac{1-(-1)^d}{2})\xi^0_C)),\eqno(4.1)$$
$$\Phi_W(\tau)=\widehat{A}(TM,\nabla^{TM})({\rm
sinh}\frac{c_0}{2})^{2n+\frac{1-(-1)^d}{2}}{\rm
ch}(\Theta_2(T_CM,(2n+\frac{1-(-1)^d}{2})\xi^0_C)),\eqno(4.2)$$
$$\Phi'_W(\tau)=\widehat{A}(TM,\nabla^{TM})({\rm
sinh}\frac{c_0}{2})^{2n+\frac{1-(-1)^d}{2}}{\rm
ch}(\Theta_3(T_CM,(2n+\frac{1-(-1)^d}{2})\xi^0_C)).\eqno(4.3)$$
Direct computations as in Section 3 and applying Chern-weil theory,
we have
$$ \Phi_L(\nabla^{TM},\nabla^{\xi^0},\tau)=\sqrt{2}^{2d-1} (\pi\sqrt{-1})^{2n+\frac{1-(-1)^d}{2}}{\rm det}^{\frac{1}{2}}
\left(\frac{R^{TM}}{4\pi^2}\frac{\theta'(0,\tau)}{\theta(\frac{R^{TM}}{4\pi^2},\tau)}\frac
{\theta_1(\frac{R^{TM}}{4\pi^2},\tau)}{\theta_1(0,\tau)}\right)$$
$$\cdot\left(\frac{\theta(u',\tau)}{\theta'(0,\tau)}\frac
{\theta_1(0,\tau)}{\theta_1(u',\tau)}\right)^{2n+\frac{1-(-1)^d}{2}}
,\eqno(4.4)$$
$$ \Phi_W(\nabla^{TM},\nabla^{\xi^0},\tau)=(\pi\sqrt{-1})^{2n+\frac{1-(-1)^d}{2}}{\rm det}^{\frac{1}{2}}
\left(\frac{R^{TM}}{4\pi^2}\frac{\theta'(0,\tau)}{\theta(\frac{R^{TM}}{4\pi^2},\tau)}\frac
{\theta_2(\frac{R^{TM}}{4\pi^2},\tau)}{\theta_2(0,\tau)}\right)$$
$$\cdot\left(\frac{\theta(u',\tau)}{\theta'(0,\tau)}\frac
{\theta_2(0,\tau)}{\theta_2(u',\tau)}\right)^{2n+\frac{1-(-1)^d}{2}}
,\eqno(4.5)$$
$$ \Phi'_W(\nabla^{TM},\nabla^{\xi^0},\tau)=(\pi\sqrt{-1})^{2n+\frac{1-(-1)^d}{2}}{\rm det}^{\frac{1}{2}}
\left(\frac{R^{TM}}{4\pi^2}\frac{\theta'(0,\tau)}{\theta(\frac{R^{TM}}{4\pi^2},\tau)}\frac
{\theta_3(\frac{R^{TM}}{4\pi^2},\tau)}{\theta_3(0,\tau)}\right)$$
$$\cdot\left(\frac{\theta(u',\tau)}{\theta'(0,\tau)}\frac
{\theta_3(0,\tau)}{\theta_3(u',\tau)}\right)^{2n+\frac{1-(-1)^d}{2}}
,\eqno(4.6)$$ where $u'=\frac{\sqrt{-1}{\rm Pf}(R^{\xi^0})}{2\pi}$
(cf. [Z]). As in [CH] and [W], we transgress $\Phi_L,\Phi_W,\Phi'_W$
and get the following forms:\\
$CS\Phi_L(\nabla_0^{TM},\nabla_1^{TM},\nabla^{\xi^0},\tau)$
$$:= \frac{\sqrt{2}}{8\pi^2}\int_0^1
\Phi_L(\nabla_t^{TM},\nabla^{\xi^0},\tau){\rm tr}\left[A \left(
\frac{1}{\frac{R_t^{TM}}{4\pi^2}}-\frac{\theta'(\frac{R_t^{TM}}{4\pi^2},\tau)}
{\theta(\frac{R_t^{TM}}{4\pi^2},\tau)}+\frac{\theta_1'(\frac{R_t^{TM}}{4\pi^2},\tau)}
{\theta_1(\frac{R_t^{TM}}{4\pi^2},\tau)}\right)\right]dt;\eqno(4.7)$$
$CS\Phi_W(\nabla_0^{TM},\nabla_1^{TM},\nabla^{\xi^0},\tau)$
$$:= \frac{1}{8\pi^2}\int_0^1
\Phi_W(\nabla_t^{TM},\nabla^{\xi^0},\tau){\rm tr}\left[A \left(
\frac{1}{\frac{R_t^{TM}}{4\pi^2}}-\frac{\theta'(\frac{R_t^{TM}}{4\pi^2},\tau)}
{\theta(\frac{R_t^{TM}}{4\pi^2},\tau)}+\frac{\theta_2'(\frac{R_t^{TM}}{4\pi^2},\tau)}
{\theta_2(\frac{R_t^{TM}}{4\pi^2},\tau)}\right)\right]dt;\eqno(4.8)$$
$CS\Phi_W'(\nabla_0^{TM},\nabla_1^{TM},\nabla^{\xi^0},\tau)$
$$:= \frac{1}{8\pi^2}\int_0^1
\Phi_W'(\nabla_t^{TM},\nabla^{\xi^0},\tau){\rm tr}\left[A \left(
\frac{1}{\frac{R_t^{TM}}{4\pi^2}}-\frac{\theta'(\frac{R_t^{TM}}{4\pi^2},\tau)}
{\theta(\frac{R_t^{TM}}{4\pi^2},\tau)}+\frac{\theta_3'(\frac{R_t^{TM}}{4\pi^2},\tau)}
{\theta_3(\frac{R_t^{TM}}{4\pi^2},\tau)}\right)\right]dt,\eqno(4.9)$$
which lie in $\Omega^{\rm odd}(M,{\bf C})[[q^{\frac{1}{2}}]]$ and
their top components represent elements in $H^{2d-1}(M,{\bf
C})[[q^{\frac{1}{2}}]]$. We have the following results.\\

\noindent {\bf Theorem 4.1} {\it Let $M$ be a $2d-1$ dimensional
manifold and  $\nabla_0^{TM},~\nabla_1^{TM}$ be two connections on
$TM$ and $\xi^0$ be a two dimensional oriented Euclidean real vector
bundle with a Euclidean connection $\nabla^{\xi^0}$, then we have\\
\noindent 1)
$\{CS\Phi_L(\nabla_0^{TM},\nabla_1^{TM},\nabla^{\xi^0},\tau)\}^{(2d-1)}$
is a modular form of weight $d-(2n+\frac{1-(-1)^d}{2})$ over $\Gamma_0(2)$;\\
$\{CS\Phi_W(\nabla_0^{TM},\nabla_1^{TM},\nabla^{\xi^0},\tau)\}^{(2d-1)}$
is a modular form of weight $d-(2n+\frac{1-(-1)^d}{2})$ over $\Gamma^0(2);$\\
$\{CS\Phi_W'(\nabla_0^{TM},\nabla_1^{TM},\nabla^{\xi^0},\tau)\}^{(2d-1)}$
is a modular form of weight $d-(2n+\frac{1-(-1)^d}{2})$ over $\Gamma_\theta(2).$\\
2) The following equalities hold,}
$$\{CS\Phi_L(\nabla_0^{TM},\nabla_1^{TM},\nabla^{\xi^0},-\frac{1}{\tau})\}^{(2d-1)}
=2^d\tau^{d-(2n+\frac{1-(-1)^d}{2})}\{CS\Phi_W(\nabla_0^{TM},\nabla_1^{TM},\nabla^{\xi^0},\tau)\}^{(2d-1)},$$
$$CS\Phi_W(\nabla_0^{TM},\nabla_1^{TM},\nabla^{\xi^0},\tau+1)
=CS\Phi_W'(\nabla_0^{TM},\nabla_1^{TM},\nabla^{\xi^0},\tau).$$

\noindent{\bf Proof.} By (2.12)-(2.17) and (4.7)-(4.9), we have
$$\{CS\Phi_L(\nabla_0^{TM},\nabla_1^{TM},\nabla^{\xi^0},-\frac{1}{\tau})\}^{(2d-1)}=
2^d\tau^{d-(2n+\frac{1-(-1)^d}{2})}\{CS\Phi_W(\nabla_0^{TM},
\nabla_1^{TM},\nabla^{\xi^0},\tau)\}^{(2d-1)},$$
$$CS\Phi_L(\nabla_0^{TM},\nabla_1^{TM},\nabla^{\xi^0},{\tau}+1)
=CS\Phi_L(\nabla_0^{TM},\nabla_1^{TM},\nabla^{\xi^0},\tau),$$
$$\{CS\Phi_W(\nabla_0^{TM},\nabla_1^{TM},\nabla^{\xi^0},-\frac{1}{\tau})\}^{(2d-1)}
=2^{-d}\tau^{d-(2n+\frac{1-(-1)^d}{2})}
\{CS\Phi_L(\nabla_0^{TM},\nabla_1^{TM},\nabla^{\xi^0},\tau)\}^{(2d-1)},$$
$$CS\Phi_W(\nabla_0^{TM},\nabla_1^{TM},\nabla^{\xi^0},{\tau}+1)
=CS\Phi_W'(\nabla_0^{TM},\nabla_1^{TM},\nabla^{\xi^0},\tau),$$
$$\{CS\Phi_W'(\nabla_0^{TM},\nabla_1^{TM},\nabla^{\xi},-\frac{1}{\tau})\}^{(2d-1)}
=(\tau)^{d-(2n+\frac{1-(-1)^d}{2})}\{CS\Phi_W'(\nabla_0^{TM},\nabla_1^{TM},\nabla^{\xi^0},\tau)\}^{(2d-1)},$$
$$CS\Phi_W'(\nabla_0^{TM},\nabla_1^{TM},\nabla^{\xi^0},{\tau}+1)
=CS\Phi_W(\nabla_0^{TM},\nabla_1^{TM},\nabla^{\xi^0},\tau).\eqno(4.10)$$
From (4.10), we can prove Theorem 4.1. $\Box$\\

 \indent Let $d=4$ and $n=1$, i.e. for $7$-dimensional manifold,
$\{CS\Phi_L(\nabla_0^{TM},\nabla_1^{TM},\nabla^{\xi^0},\tau)\}^{7}$
is a modular form of weight $2$ over $\Gamma_0(2)$. Set
$A=\nabla_1^{TM}-\nabla_0^{TM}.$ Using similar discussions in [CH,
p.15], we get
$$CS\Phi_L(\nabla_0^{TM},\nabla_1^{TM},\nabla^{\xi^0},\tau) =-\frac
{\delta_1}{6\pi^2}c^2_0{\rm
tr}\left[A[\nabla_0^{TM},\nabla_1^{TM}]+\frac{2}{3}A\wedge A\wedge
A\right].\eqno(4.11)$$ Similarly, we obtain that
$$CS\Phi_W(\nabla_0^{TM},\nabla_1^{TM},\nabla^{\xi^0},\tau) =-\frac
{\delta_2}{96\pi^2}c^2_0{\rm
tr}\left[A[\nabla_0^{TM},\nabla_1^{TM}]+\frac{2}{3}A\wedge A\wedge
A\right].\eqno(4.12)$$
$$CS\Phi_W'(\nabla_0^{TM},\nabla_1^{TM},\nabla^{\xi^0},\tau) =-\frac
{\delta_3}{96\pi^2}c^2_0{\rm
tr}\left[A[\nabla_0^{TM},\nabla_1^{TM}]+\frac{2}{3}A\wedge A\wedge
A\right].\eqno(4.12)$$

From Theorem 4.1, we can imply some twisted cancellation formulas
for odd dimensional manifolds. For example, let $d=6$ and $n=1$,
i.e.  $M$ be $11$ dimensional. We have that
$\{CS\Phi_L(\nabla_0^{TM},\nabla_1^{TM},\nabla^{\xi^0},\tau)\}^{(11)}$
is a modular form of weight $4$ over $\Gamma_0(2),$
$\{CS\Phi_W(\nabla_0^{TM},\nabla_1^{TM},\nabla^{\xi^0},\tau)\}^{(11)}$
is a modular form of weight $4$ over $\Gamma^0(2)$ and
$$\{CS\Phi_L(\nabla_0^{TM},\nabla_1^{TM},\nabla^{\xi^0},-\frac{1}{\tau})\}^{(11)}
=2^6\tau^{4}\{CS\Phi_W(\nabla_0^{TM},\nabla_1^{TM},\nabla^{\xi^0},\tau)\}^{(11)}.\eqno(4.13)$$
By Lemma 2.2, we have
$$\{CS\Phi_W(\nabla_0^{TM},\nabla_1^{TM},\nabla^{\xi^0},\tau)\}^{(11)}
=z_0(8\delta_2)^2+z_1\varepsilon_2,\eqno(4.14))$$ and by (2.19) and
Theorem (4.13),
$$\{CS\Phi_L(\nabla_0^{TM},\nabla_1^{TM},\nabla^{\xi},\tau)\}^{(11)}
=2^6[z_0(8\delta_1)^2+z_1\varepsilon_1].\eqno(4.15)$$ By comparing
the $q^{\frac{1}{2}}$-expansion coefficients in (4.14), we get

$$z_0=\left\{\int_0^1\widehat{A}(TM,\nabla^{TM}_t)({\rm
sinh}\frac{c_0}{2})^2{\rm
tr}\left[A\left(\frac{1}{2R_t^{TM}}-\frac{1}{8\pi{\rm
tan}{\frac{R^{TM}_t}{4\pi}}}\right)\right]dt\right\}^{(11)},\eqno(4.16)$$
$$z_1=\left\{\int_0^1\widehat{A}(TM,\nabla^{TM}_t)({\rm
sinh}\frac{c_0}{2})^2\left(-{\rm
ch}(T_CM,\nabla^{T_CM}_t)+2(e^{c_0}+e^{-c_0})-41\right)\right.$$
$$\left.\times{\rm tr}\left[A\left(\frac{1}{2R_t^{TM}}-\frac{1}{8\pi{\rm
tan}{\frac{R^{TM}_t}{4\pi}}}\right)\right]dt\right.+$$
$$\left.\int_0^1\widehat{A}(TM,\nabla^{TM}_t)({\rm
sinh}\frac{c_0}{2})^2{\rm tr}\left[\frac{A}{2\pi}{\rm
sin}\frac{R^{TM}_t}{4\pi}\right]dt\right\}^{(11)}.\eqno(4.17)$$
Plugging (4.16) and (4.17) into (4.15) and comparing the constant
terms of both sides, we obtain that
$$\left\{\int_0^1\sqrt{2}\widehat{L}(TM,\nabla^{TM}_t)\frac{{\rm
sinh}^2\frac{c_0}{2}}{{\rm cosh}^2\frac{c_0}{2}}{\rm
tr}\left[A\left(\frac{1}{2R_t^{TM}}-\frac{1}{4\pi{\rm
sin}\frac{R^{TM}_t}{2\pi}}\right)\right]\right\}^{(11)}=2^2(2^6z_0+z_1),\eqno(4.18)$$
so we have the following $11$-dimensional analogue of the twisted miraculous cancellation formula.\\

\noindent {\bf Corollary 3.3} {\it The following equality holds}
$$\left\{\int_0^1\sqrt{2}\widehat{L}(TM,\nabla^{TM}_t)\frac{{\rm
sinh}^2{\frac{c_0}{2}}}{{\rm cosh}^2{\frac{c_0}{2}}}{\rm
tr}\left[A\left(\frac{1}{2R_t^{TM}}-\frac{1}{4\pi{\rm
sin}{\frac{R^{TM}_t}{2\pi}}}\right)\right]\right\}^{(11)}$$
$$=4\left\{\int_0^1\widehat{A}(TM,\nabla^{TM}_t)({\rm
sinh}\frac{c_0}{2})^2\left(-{\rm
ch}(T_CM,\nabla^{T_CM}_t)+2(e^{c_0}+e^{-c_0})+23\right)\right.$$
$$\left.\times{\rm tr}\left[A\left(\frac{1}{2R_t^{TM}}-\frac{1}{8\pi{\rm
tan}{\frac{R^{TM}_t}{4\pi}}}\right)\right]dt\right.+$$
$$\left.\int_0^1\widehat{A}(TM,\nabla^{TM}_t)({\rm
sinh}\frac{c_0}{2})^2{\rm tr}\left[\frac{A}{2\pi}{\rm
sin}\frac{R^{TM}_t}{4\pi}\right]dt\right\}^{(11)}.\eqno(4.19)$$

let $d=5$ and $n=0$, i.e.  $M$ be $9$ dimensional. Using similar
discussions, we have\\

\noindent {\bf Corollary 3.4} {\it The following equality holds}
$$\left\{\int_0^1\sqrt{2}\widehat{L}(TM,\nabla^{TM}_t)\frac{{\rm
sinh}{\frac{c_0}{2}}}{{\rm cosh}{\frac{c_0}{2}}}{\rm
tr}\left[A\left(\frac{1}{2R_t^{TM}}-\frac{1}{4\pi{\rm
sin}{\frac{R^{TM}_t}{2\pi}}}\right)\right]\right\}^{(9)}$$
$$=2\left\{\int_0^1\widehat{A}(TM,\nabla^{TM}_t)({\rm
sinh}\frac{c_0}{2})\left(-{\rm
ch}(T_CM,\nabla^{T_CM}_t)+(e^{c_0}+e^{-c_0})+23\right)\right.$$
$$\left.\times{\rm tr}\left[A\left(\frac{1}{2R_t^{TM}}-\frac{1}{8\pi{\rm
tan}{\frac{R^{TM}_t}{4\pi}}}\right)\right]dt\right.+$$
$$\left.\int_0^1\widehat{A}(TM,\nabla^{TM}_t)({\rm
sinh}\frac{c_0}{2}){\rm tr}\left[\frac{A}{2\pi}{\rm
sin}\frac{R^{TM}_t}{4\pi}\right]dt\right\}^{(9)}.\eqno(4.20)$$
\\

\noindent{\bf Acknowledgement}~This work was supported by NSFC No.
10801027.\\

\noindent {\bf References}\\

\noindent [AW] L. Alvarez-Gaum\'{e}, E. Witten, Graviational
anomalies,
{\it Nucl. Phys.} B234 (1983), 269-330.\\
\noindent [Ch] K. Chandrasekharan, {\it Elliptic Functions},
Spinger-Verlag, 1985. \\
\noindent [CH] Q. Chen, F. Han, Elliptic genera, transgression and
loop space Chern-Simons form, arXiv:0605366.\\
\noindent [HZ1]F. Han, W. Zhang, ${\rm Spin}^c$-manifold and
elliptic genera,
{\it C. R. Acad. Sci. Paris Serie I.,} 336 (2003), 1011-1014.\\
\noindent [HZ2] F. Han, W. Zhang, Modular invariance, characteristic
numbers and eta Invariants,  {\it J.
Diff. Geom.} 67 (2004), 257-288.\\
 \noindent [HH] F. Han, X. Huang, Even dimensional manifolds and generalized anomaly cancellation
formulas , {\it Trans. AMS }359 (2007), No. 11, 5365-5382.\\
\noindent [Li] K. Liu, Modular invariance and characteristic
numbers. {\it Commu.Math. Phys.} 174 (1995), 29-42.\\
\noindent [W] Y. Wang, Transgression and twisted anomaly
cancellation formulas on odd dimensional manifolds, preprint.\\
 \noindent [Z] W. Zhang, {\it Lectures on Chern-weil Theory
and Witten Deformations.} Nankai Tracks in Mathematics Vol. 4, World
Scientific, Singapore, 2001.\\

\indent School of Mathematics and Statistics , Northeast
Normal University, Changchun, Jilin 130024, China ;\\

 \indent E-mail: {\it wangy581@nenu.edu.cn}
\end {document}